\documentclass[12pt]{amsart}
\usepackage{amssymb}
\usepackage{amsfonts}
\usepackage{latexsym}
\usepackage{amscd}
\usepackage[mathscr]{euscript}
\usepackage{xy} \xyoption{all}

\newcommand{\Z}{{\mathbb{Z}}}

\newcommand{\C}{{\mathbb{C}}}

\newcommand{\ol}{\overline}
\newcommand{\uloopr}[1]{\ar@'{@+{[0,0]+(-4,5)}@+{[0,0]+(0,10)}@+{[0,0] +(4,5)}}^{#1}}
\newcommand{\uloopd}[1]{\ar@'{@+{[0,0]+(5,4)}@+{[0,0]+(10,0)}@+{[0,0]+ (5,-4)}}^{#1}}
\newcommand{\dloopr}[1]{\ar@'{@+{[0,0]+(-4,-5)}@+{[0,0]+(0,-10)}@+{[0, 0]+(4,-5)}}_{#1}}
\newcommand{\dloopd}[1]{\ar@'{@+{[0,0]+(-5,4)}@+{[0,0]+(-10,0)}@+{[0,0 ]+(-5,-4)}}_{#1}}

\newcommand{\luloop}[1]{\ar@'{@+{[0,0]+(-8,2)}@+{[0,0]+(-10,10)}@+{[0, 0]+(2,2)}}^{#1}}

\newtheorem{lemma}{Lemma}[section]
\newtheorem{corollary}[lemma]{Corollary}
\newtheorem{theorem}[lemma]{Theorem}
\newtheorem{proposition}[lemma]{Proposition}
\newtheorem{remark}[lemma]{Remark}
\newtheorem{definition}[lemma]{Definition}
\newtheorem{example}[lemma]{Example}

\begin{document}

\title[Nonstable $K$-theory for graph algebras]{Nonstable $K$-theory for graph algebras}%
\author{P. Ara}
\address{Departament de Matem\`atiques, Universitat Aut\`onoma de Barcelona,
08193 Bellaterra (Barcelona), Spain.} \email{para@mat.uab.cat}
\author{M.A. Moreno}\author{E. Pardo}
\address{Departamento de Matem\'aticas, Universidad de C\'adiz,
Apartado 40, 11510 Puerto Real (C\'adiz),
Spain.}\email{mariangeles.moreno@uca.es}\email{enrique.pardo@uca.es}
\urladdr{http://www.uca.es/dept/matematicas/PPersonales/PardoEspino/EMAIN.HTML}
\thanks{The first author was partially supported by the DGI
and European Regional Development Fund, jointly, through Project
BFM2002-01390, the second and the third  by the DGI and European
Regional Development Fund, jointly, through Project MTM2004-00149
and by PAI III grant FQM-298 of the Junta de Andaluc\'{\i}a. Also,
the first and third authors are partially supported by the
Comissionat per Universitats i Recerca de la Generalitat de
Catalunya.} 
\date{\today}

\begin{abstract}
We compute the monoid $V(L_K(E))$ of isomorphism classes of
finitely generated projective modules over certain graph algebras
$L_K(E)$, and we show that this monoid satisfies the refinement
property and separative cancellation. We also show that there is a
natural isomorphism between the lattice of graded ideals of
$L_K(E)$ and the lattice of order-ideals of $V(L_K(E))$. When $K$
is the field $\mathbb C$ of complex numbers, the algebra
$L_{\mathbb C}(E)$ is a dense subalgebra of the graph
$C^*$-algebra $C^*(E)$, and we show that the inclusion map induces
an isomorphism between the corresponding monoids. As a
consequence, the graph C*-algebra of any row-finite graph turns
out to satisfy the stable weak cancellation property.
\end{abstract}
\maketitle

\section{Introduction}

The Cuntz-Krieger algebras $\mathcal O _A$, introduced by Cuntz
and Krieger \cite{CK} in 1980, constitute a prominent class of
$C^*$-algebras. The algebras $\mathcal O _A$ were originally
associated to a finite matrix $A$ with entries in $\{ 0,1\}$, but
it was quickly realized that they could also be viewed as the
$C^*$-algebras of a finite directed graph \cite{Wat}. These
$C^*$-algebras, as well as those arising from various infinite
graphs, have been the subject of much investigation (see e.g.
\cite{BPRS}, \cite{DHS}, \cite{KPR}, \cite{JP}, \cite{RS}).
Although Raeburn and Szyma\'nski \cite{RS} have computed the $K_0$
and $K_1$-groups of a graph $C^*$-algebra $C^*(E)$ associated with
any row-finite graph $E$, the actual structure of the monoid
$V(C^*(E))$ of Murray-von Neumann equivalence classes of
projections in matrix algebras over $C^*(E)$ seems to remain
unnoticed. One of the goals of this paper is to fill this gap.
Another major goal is to show some nice decomposition and
cancellation properties of projections over graph $C^*$-algebras,
which also hold for purely algebraic versions of them.

Any graph $C^*$-algebra $C^*(E)$ is the completion, in an
appropriate norm, of a certain $*$-subalgebra $L_{\mathbb C}(E)$,
which is just the $*$-subalgebra generated by the canonical
projections and partial isometries that generate $C^*(E)$ as a
$C^*$-algebra. We show that the natural inclusion $\psi \colon
L_{\mathbb C}(E)\to C^*(E)$ induces a monoid isomorphism
$V(\psi)\colon V(L_{\mathbb C}(E))\to V(C^*(E))$ (Theorem
\ref{C*algebras}). In this algebraic vein, similar algebras
$L_K(E)$ can be constructed over an arbitrary field $K$, and we
show that the monoid $V(L_K(E))$ does not depend on the field $K$.
The algebras $L_K(E)$ have been already considered recently by
Abrams and Aranda Pino in \cite{AA}, under the name of Leavitt
path algebras. They provide a generalization of Leavitt algebras
of type $(1,n)$, introduced by Leavitt in 1962 \cite{Leavitt},
just in the same way as graph $C^*$-algebras $C^*(E)$ provide a
generalization of Cuntz algebras.

The decomposition properties of projections in (matrix algebras
over) a $C^*$-algebra $A$ are faithfully reflected in the
structure of the monoid $V(A)$. This is an essential ingredient in
the so-called nonstable $K$-theory for $C^*$-algebras; cf.
\cite{RatC*}. A similar statement holds true in Ring Theory, where
the monoid $V(R)$ is usually described in terms of the finitely
generated projective $R$-modules; see for example \cite{AGOP} and
\cite{AF}. For a $C^*$-algebra $A$, the two versions of $V(A)$,
obtained by viewing $A$ as a $C^*$-algebra or viewing $A$ as a
plain ring, agree, see Section 2. Moreover, important information
about the lattice of ideals of a ring $R$ is faithfully codified
in the monoid $V(R)$; see for example \cite[Theorem 2.1]{FHK}. The
subsets of $V(R)$ corresponding to ideals in $R$ are the so-called
order-ideals of $V(R)$, which are the submonoids $S$ of $V(R)$
such that, for $x,y\in V(R)$, we have $x+y\in S$ if and only if
$x\in S$ and $y\in S$. Then, \cite[Theorem 2.1]{FHK} asserts that
the lattice of all order-ideals of $V(R)$ is isomorphic with the
lattice of all {\it trace ideals} of $R$.

We consider an abelian monoid $M_E$ associated with a directed
row-finite graph $E$, and we prove that this monoid is naturally
isomorphic with the monoid of isomorphism classes of finitely
generated projective modules over $L_K(E)$ (see Theorem
\ref{VLE}). This uses the nice machinery developed by Bergman in
\cite{Bergman} to compute the monoids $V(R)$ of algebras $R$
obtained by means of some universal constructions. We also show
that the monoid $M_E$ is naturally isomorphic to $V(C^*(E))$
(Theorem \ref{C*algebras}), and indeed that the natural map
$L_{\mathbb C}(E)\to C^*(E)$ induces an isomorphism $V(L_{\mathbb
C}(E))\to V(C^*(E))$, but our proof of this fact is quite
involved, basically because we do not have at our disposal a
$C^*$-version of Bergman's machinery. Rather, our proof uses the
computation in \cite{RS} of $K_0(C^*(E))$, which implies that
$K_0(C^*(E))$ agrees with the universal group of $M_E$. We then
use techniques from nonstable $K$-theory to deduce the equality of
the monoids $V(C^*(E))$ and $M_E$. As a consequence of this fact
and of our monoid theoretic study of the monoid $M_E$, we get that
$C^*(E)$ always has stable weak cancellation (Corollary
\ref{stableweak}) (equivalently, $C^*(E)$ is separative, see
Proposition \ref{swc}). An analogous result holds for all graph
algebras $L_K(E)$ (Corollary \ref{separativerings}). We remark
that various stability results for wide classes of rings and
$C^*$-algebras can be proved under the additional hypothesis of
separativity; see for example \cite{AGOP}, \cite{AGOR},
\cite{Perera}, \cite{BPns}.

The most important tools from nonstable $K$-theory we use are the
concepts of refinement and separative cancellation. These
properties are faithfully reflected in monoid theoretic properties
of the associated monoid $V(A)$. A substantial part of this paper
is devoted to establish these properties for the monoid $M_E$,
using just monoid theoretic techniques. Both concepts were defined
and studied in \cite{AGOP}. The definitions will be recalled in
Section 2.

We now summarize the contents of the rest of sections of the
paper. Section 3 contains the definition of the (Leavitt) graph
algebras $L_K(E)$ and of the monoid $M_E$ associated with a
row-finite graph $E$. The monoid $M_E$ is isomorphic to
$F_E/{\sim}$, where $F_E$ is the free abelian monoid on $E^0$ and
$\sim$ is a certain congruence on $F_E$. Our basic tool for the
monoid theoretic study of $M_E$ is a precise description of this
congruence, which is given in Section 4, which also contains the
proof of the refinement property of $M_E$. In Section 5, we shall
establish an isomorphism between the lattice $\mathcal H$ of
saturated hereditary subsets of $E^0$, the lattice of order-ideals
of $M_E$, and the lattice of graded ideals of $L_K(E)$. This
result parallels \cite[Theorem 4.1]{BPRS}, where an isomorphism
between the lattice of saturated hereditary subsets of $E^0$ and
the lattice of closed gauge-invariant ideals of the $C^*$-algebra
$C^*(E)$ is obtained. The separativity property of the monoid
$M_E$ is obtained in Section 6. Finally we show in Section 7 that
$V(C^*(E))$ is naturally isomorphic with $M_E$. This result,
together with all properties we have obtained for $M_E$, enables
us to conclude that $C^*(E)$ has stable weak cancellation.

\section{Basic concepts}
Our references for $K$-theory for $C^*$-algebras are \cite{Black}
and \cite{rordam}. For algebraic $K$-theory, we refer the reader
to \cite{Ros}. For a unital ring $R$, let $M_{\infty}(R)$ be the
directed union of $M_n(R)$ ($n\in\mathbb N$), where the transition
maps $M_n(R)\to
M_{n+1}(R)$ are given by $x\mapsto \left( \smallmatrix x&0\\
0&0\endsmallmatrix \right)$. We define $V(R)$ to be the set of
isomorphism classes (denoted $[P]$) of finitely generated
projective left $R$-modules, and we endow $V(R)$ with the
structure of a commutative monoid by imposing the operation
$$[P]+ [Q] := [P\oplus Q]$$ for any isomorphism classes $[P]$ and
$[Q]$. Equivalently \cite[Chapter 3]{Black}, $V(R)$ can be viewed
as the set of equivalence classes $V(e)$ of idempotents $e$ in
$M_\infty(R)$ with the operation
$$V(e)+V(f) :=
V\bigl(  \left( \smallmatrix e&0\\ 0&f
\endsmallmatrix \right) \bigr)$$
for idempotents $e,f\in M_\infty(R)$. The group $K_0(R)$ of a
unital ring $R$ is the universal group of $V(R)$. Recall that, as
any universal group of an abelian monoid, the group $K_0(R)$ has a
standard structure of partially pre-ordered abelian group. The set
of positive elements in $K_0(R)$ is the image of $V(R)$ under the
natural monoid homomorphism $V(R)\to K_0(R)$. Whenever $A$ is a
$C^*$-algebra, the monoid $V(A)$ agrees with the monoid of
Murray-von Neumann equivalence classes of projections in
$M_{\infty}(A)$; see \cite[4.6.2 and 4.6.4]{Black} or
\cite[Exercise 3.11]{rordam}. It follows that the algebraic
version of $K_0(A)$ coincides with the operator-theoretic one.
\medskip

We now review some important decomposition and cancellation
properties concerning finitely generated projective modules. In
the context of $C^*$-algebras, these are equivalent to
corresponding statements for projections, as in \cite[Section
7]{AGOP}.

Let $FP(R)$ be the class of finitely generated projective modules
over a ring $R$. We say that $FP(R)$ satisfies the {\em refinement
property} if whenever $A_1,A_2,B_1,B_2\in FP(R)$ satisfy
$A_1\oplus A_2\cong B_1\oplus B_2$, there exist decompositions
$A_i=A_{i1}\oplus A_{i2}$ for $i=1,2$ such that $A_{1j}\oplus
A_{2j}\cong B_j$ for $j=1,2$.

It was proved in \cite[Proposition 1.2]{AGOP} that every exchange
ring satisfies the refinement property. Among $C^*$-algebras, it
is worth to mention that every $C^*$-algebra with real rank zero
(\cite{BP}) satisfies the refinement property. This is a theorem
of Zhang \cite[Theorem 3.2]{zhang}. It can also be seen as a
consequence of the above mentioned result on exchange rings, since
every $C^*$-algebra of real rank zero is an exchange ring
\cite[Theorem 7.2]{AGOP}.

An abelian monoid $M$ is a {\em refinement monoid} if whenever
$a+b=c+d$ in $M$, there exist $x,y,z,t\in M$ such that $a=x+y$ and
$b=z+t$ while $c=x+z$ and $d=y+t$. It is clear that $V(R)$ is a
refinement monoid if and only if the class $FP(R)$ satisfies the
refinement property. We will show that this is the case when
$R=L_K(E)$ or $R=C^*(E)$.
\smallskip

Now we discuss the concept of separative cancellation. We say that
a ring $R$ is {\em separative} in case it satisfies the following
property: If $A,B,C\in FP(R)$ satisfy $A\oplus C\cong B\oplus C$
and $C$ is isomorphic to direct summands of both $nA$ and $nB$ for
some $n\in \mathbb N$, then $A\cong B$.

Many rings are separative. Indeed it is an outstanding open
question to determine whether all exchange rings are separative.
In the context of $C^*$-algebras, it is not known whether all
$C^*$-algebras of real rank zero are separative. We will show that
all graph $C^*$-algebras $C^*(E)$ and all Leavitt graph algebras
$L_K(E)$ are separative.

This concept is closely related to the concept of weak
cancellation, introduced by L.G. Brown in \cite{brown}. See also
\cite{BPns}, where many extremally rich $C^*$-algebras are shown
to have weak cancellation. Following \cite{brown} and \cite{BPns},
we say that a $C^*$-algebra $A$ has {\em weak cancellation} if any
pair of projections $p,q$ in $A$ that generate the same closed
ideal $I$ in $A$ and have the same image in $K_0(I)$ must be
Murray-von Neumann equivalent in $A$ (hence in $I$). If $M_n(A)$
has weak cancellation for every $n$, then we say that $A$ has {\em
stable weak cancellation}. It is an open problem whether every
extremally rich $C^*$-algebra satisfies weak cancellation. By
\cite[Theorem 2.11]{BPns}, every extremally rich $C^*$-algebra of
real rank zero has stable weak cancellation.

If $P$ and $Q$ are projections in $M_{\infty}(A)$, we will use the
symbol $P\sim Q$ to indicate that they are (Murray-von Neumann)
equivalent, that is, there is a partial isometry $W$ in
$M_{\infty}(A)$ such that $W^*W=P$ and $WW^*=Q$. Similarly, we
will write $P\lesssim Q$ in case $P$ is equivalent to a projection
$Q'$ such that $Q'=Q'Q$. We will write $P\oplus Q$ for the
block-diagonal matrix $\text{diag}(P,Q)$, and we will denote by
$n\cdot P$ the direct sum of $n$ copies of $P$.

\begin{proposition}
\label{swc} Let $A$ be a $C^*$-algebra. Then $A$ has stable weak
cancellation if and only if $A$ is separative.
\end{proposition}

\begin{proof}
The proof is straightforward, taking into account the following
fact: Two projections $P, Q\in M_{\infty}(A)$ whose respective
entries generate the same closed ideal $I$ of $A$ have the same
image in $K_0(I)$ if and only if there is a projection $E\in
M_{\infty}(I)$ such that $P\oplus E\sim Q\oplus E$. (Note that a
projection  $E$ belongs to $M_{\infty}(I)$ if and only if
$E\lesssim n\cdot P $ and $E\lesssim n\cdot Q$ for some $n\ge 1$.)
\end{proof}

There is a canonical pre-order on any abelian monoid $M$, which is
sometimes called the algebraic pre-order of $M$. This pre-order is
defined by setting $x\le y$ if and only if there is $z\in M$ such
that $y=x+z$. This is the only pre-order that we will consider in
this paper for a monoid.

An abelian monoid $M$ is said to be {\em separative} \cite{AGOP}
in case $M$ satisfies the following condition: If $a,b,c\in M$
satisfy $a+c=b+c$ and $c\le na$ and $c\le nb$ for some $n\in
\mathbb N$, then $a=b$. It is clear that a ring $R$ is separative
if and only if the monoid $V(R)$ is separative.

\section{Graph algebras and graph monoids}
A directed graph $E$ consists of a vertex set $E^0$, an edge set
$E^1$, and maps $r,s: E^1\longrightarrow E^0$ describing the range
and source of edges. We say that $E$ is a {\em row-finite graph}
if each row in its adjacency matrix $A_E=(A(v,w))_{v,w\in E^0}$
has only a finite number of nonzero entries, where $A(v,w)$ is the
number of edges going from $v$ to $w$. This amounts to saying that
each vertex in $E$ emits only a finite number of edges.

Let $E=(E^0,E^1)$ be a row-finite graph, and let $K$ be a field.
We define the {\em graph $K$-algebra} $L_K(E)$ associated with $E$
as the $K$-algebra generated by a set $\{p_v\mid v\in E^0\}$
together with a set $\{x_e,y_e\mid e\in E^1\}$, which satisfy the
following relations:

(1) $p_vp_{v'}=\delta _{v,v'}p_v$ for all $v,v'\in E^0$.

(2) $p_{s(e)}x_e=x_ep_{r(e)}=x_e$ for all $e\in E^1$.

(3) $p_{r(e)}y_e=y_ep_{s(e)}=y_e$ for all $e\in E^1$.

(4) $y_ex_{e'}=\delta _{e,e'}p_{r(e)}$ for all $e,e'\in E^1$.

(5) $p_v=\sum _{\{ e\in E^1\mid s(e)=v \}}x_ey_e$ for every $v\in
E^0$ that emits edges.

Observe that relation (1) says that $\{ p_v\mid v\in E^0\}$ is a
set of pairwise orthogonal idempotents. Note also that the above
relations imply that $\{x_ey_e\mid e\in E^1\}$ is a set of
pairwise orthogonal idempotents in $L_K(E)$. If $E$ is a finite
graph then we have $\sum _{v\in E^0} p_v=1$. In general the
algebra $L_K(E)$ is not unital, but it can be written as a direct
limit of unital graph algebras (with non-unital transition maps),
so that it is an algebra with local units. To show this, we first
observe the {\em functoriality property} of the construction, as
follows. Recall that a graph homomorphism $f\colon E=(E^0,E^1)\to
F=(F^0,F^1)$ is given by two maps $f^0\colon E^0\to F^0$ and
$f^1\colon E^1\to F^1$ such that $r_F(f^1(e))=f^0(r_E(e))$ and
$s_F(f^1(e))=f^0(s_E(e))$ for every $e\in E^1$. We say that a
graph homomorphism $f$ is {\em complete} in case $f^0$ is
injective and $f^1$ restricts to a bijection from $s_E^{-1}(v)$
onto $s_F^{-1}(f^0(v))$ for every $v\in E^0$ such that $v$ emits
edges. Note that under the above assumptions, the map $f^1$ must
also be injective. Let us consider the category $\mathcal G$ whose
objects are all the row-finite graphs and whose morphisms are the
complete graph homomorphisms. It is easy to check that the
category $\mathcal G$ admits direct limits. If $\{X_i\}_{i\in I}$
is a directed system in the category $\mathcal G$ and
$X=\varinjlim _{i\in I}X_i$, let us denote by $\psi _i\colon
X_i\to X$ the canonical direct limit homomorphisms. Then the
graphs $\psi _i(X_i)$ are complete subgraphs of $X$, $\psi
_i(X_i)$ is a complete subgraph of $\psi _j(X_j)$ whenever $i\le
j$, and $X$ is the union of the family of subgraphs $(\psi
_i(X_i))_{i\in I}$ (that is, $X^0=\cup _{i\in I}\psi _i^0(X^0_i)$
and $X^1=\cup _{i\in I}\psi _i^1(X^1_i)$).

In order to simplify notation, the $K$-algebra $L_K(E)$ will be
sometimes denoted by $L(E)$.

\begin{lemma}
\label{injlim} Every row-finite graph $E$ is a direct limit in the
category $\mathcal G$ of a directed system of finite graphs.
\end{lemma}

\begin{proof}
Clearly, $E$ is the union of its finite subgraphs. Let $X$ be a
finite subgraph of $E$. Define a finite subgraph $Y$ of $E$ as
follows:
$$Y^0=X^0\cup \{r_E(e)\mid e\in E^1 \text{ and } s_E(e)\in X^0\}$$
and
$$Y^1=\{ e\in E^1\mid s_E(e)\in X^0\}.$$ Then the vertices
of $Y$ that emit edges are exactly the vertices of $X$ that emit
edges in $E$, and if $v$ is one of these vertices, then
$s_E^{-1}(v)=s_{Y}^{-1}(v)$. This shows that the map $Y\to E$ is a
complete graph homomorphism, and clearly $X\subseteq Y$. If $Y_1$
and $Y_2$ are two complete subgraphs of $E$ and $Y_1$ is a
subgraph of $Y_2$, then the inclusion map $Y_1\to Y_2$ is clearly
a complete graph homomorphism.

Since the union of a finite number of finite complete subgraphs of
$E$ is again a finite complete subgraph of $E$, it follows that
$E$ is the direct limit in the category $\mathcal G$ of the
directed family of its finite complete subgraphs.
\end{proof}

\begin{lemma}
\label{algdirectlimit} \label{mondirectlimit}  The assignment
$E\mapsto L_K(E)$ can be extended to a functor $L_K$ from the
category $\mathcal G$ of row-finite graphs and complete graph
homomorphisms to the category of $K$-algebras and (not necessarily
unital) algebra homomorphisms. The functor $L_K$ is continuous,
that is, it commutes with direct limits. It follows that every
graph algebra $L_K(E)$ is the direct limit of graph algebras
corresponding to finite graphs.
\end{lemma}

\begin{proof}
If $f\colon E\to F$ is a complete graph homomorphism, then $f$
induces an algebra homomorphism $L(f)\colon L_K(E)\to L_K(F)$, as
follows. Set $L(f)(p_v)=p_{f^0(v)}$ and $L(f)(x_e)=x_{f^1(e)}$ and
$L(f)(y_e)=y_{f^1(e)}$ for $v\in E^0$ and $e\in E^1$. Since $f^0$
is injective, relation (1) is preserved under $L(f)$. Relations
(2), (3) are clearly preserved, relation (4) is preserved because
$f^1$ is injective, and relation (5) is preserved because $f^1$
restricts to a bijection from $s_E^{-1}(v)$ onto
$s_F^{-1}(f^0(v))$ for every $v\in E^0$ such that $v$ emits edges.

The algebra $L_K(E)$ is the algebra generated by a universal
family of elements $\{p_v,x_e,y_e\mid v\in E^0,e\in E^1\}$
satisfying relations (1)--(5). If $X=\varinjlim _{i\in I}X_i$ in
the category $\mathcal G$, then, as observed above, we can think
that $\{X_i\}_{i\in I}$ is a directed family of complete subgraphs
of $X$, and the union of the graphs $X_i$ is $X$. For a
$K$-algebra $A$, a compatible set of $K$-algebra homomorphisms
$L_K(X_i)\to A$, $i\in I$, determines, and is determined by, a set
of elements $\{p'_v,x'_e,y'_e\mid v\in E^0,e\in E^1\}$ in $A$
satisfying conditions (1)--(5). It follows that $L_K(E)=\varinjlim
_{i\in I}L_K(X_i)$, as desired. The last statement follows now
from Lemma \ref{injlim}.
\end{proof}

The graph $C^*$-algebra $C^*(E)$ is the $C^*$-algebra generated by
a universal Cuntz-Krieger $E$-family $\{ P_v,S_e \mid v\in E^0,
e\in E^1 \}$, see \cite[Theorem 1.2]{KPR}. By definition, a
Cuntz-Krieger $E$-family in a $C^*$-algebra $A$ consists of a set
$\{ P_v\mid v\in E^0\}$ of pairwise orthogonal projections in $A$
and a set $\{ S_e\mid e\in E^1\}$ of partial isometries in $A$
such that
$$ S^*_eS_e=P_{r(e)}\quad \text{ for }e\in E^1,\qquad \text{and}\qquad
P_v=\sum _{\{ e\in E^1\mid s(e)=v\}}S_eS_e^* \quad \text{ for
}v\in E^0 .
$$
Therefore the same proof as in Lemma \ref{algdirectlimit} can be
applied to the case of $C^*$-algebras:

\begin{lemma}
\label{C*algdirectlimit} \label{mondirectlimit}  The assignment
$E\mapsto C^*(E)$ can be extended to a continuous functor from the
category $\mathcal G$ of row-finite graphs and complete graph
homomorphisms to the category of $C^*$-algebras and
$*$-homomorphisms. Every graph $C^*$-algebra $C^*(E)$ is the
direct limit of graph $C^*$-algebras associated with finite
graphs.\qed
\end{lemma}

Now we want to compute the monoid $V(L_K(E))$ associated with the
finitely generated projective modules over the graph algebra
$L_K(E)$. Though $L_K(E)$ is not in general a unital algebra,
there is a well-defined monoid $V(L_K(E))$ associated with the
finitely generated projective left modules over $L_K(E)$. We
recall the general definition here.

Let $I$ be a non-unital $K$-algebra, and consider any unital
$K$-algebra $R$ containing $I$ as a two-sided ideal. We consider
the class $FP(I,R)$ of finitely generated projective left
$R$-modules $P$ such that $P=IP$. Then $V(I)$ is defined as the
monoid of isomorphism classes of objects in $FP(I,R)$, and does
not depend on the particular unital ring $R$ in which $I$ sits as
a two-sided ideal, as can be seen from the following alternative
description: $V(I)$ is the set of equivalence classes of
idempotents in $M_{\infty}(I)$, where $e\sim f$ in $M_{\infty}(I)$
if and only if there are $x,y\in M_{\infty}(I)$ such that $e=xy$
and $f=yx$. See \cite[page 296]{MM}.

The assignment $I\mapsto V(I)$ gives a functor from the category
of non-unital rings to the category of abelian monoids, that
commutes with direct limits. Moreover, being $L_K(E)$ a ring with
local units, it is well-known that $K_0(L_K(E))$, the $K_0$-group
of the non-unital ring $L_K(E)$, is just the enveloping group of
$V(L_K(E))$; see \cite[Proposition 0.1]{MM}.

Let $M_E$ be the abelian monoid given by the generators $\{
a_v\mid v\in E^0\}$, with the relations:
\begin{equation}
\label{(M)}\tag{M} a_v= \sum _{\{e\in E^1\mid
s(e)=v\}}a_{r(e)}\qquad \text{for every }v\in E^0 \text{ that
emits edges} .
\end{equation}

\begin{lemma}
\label{mondirectlimit}  The assignment $E\mapsto M_E$ can be
extended to a continuous functor from the category $\mathcal G$ of
row-finite graphs and complete graph homomorphisms to the category
of abelian monoids. It follows that every graph monoid $M_E$ is
the direct limit of graph monoids corresponding to finite graphs.
\end{lemma}

\begin{proof}
Every complete graph homomorphism $f\colon E\to F$ induces a
natural monoid homomorphism
$$M(f)\colon M_E\to M_F ,$$ and so we get a functor $M$ from the
category $\mathcal G$ to the category of abelian monoids. The fact
that $M$ commutes with direct limits is proven in the same way as
in Lemma \ref{algdirectlimit}.
\end{proof}

\begin{theorem}
\label{VLE} Let $E$ be a row-finite graph. Then there is a natural
monoid isomorphism $V(L_K(E))\cong M_E$. Moreover, if $E$ is
finite, then the global dimension of $L_K(E)$ is $\le 1$.
\end{theorem}

\begin{proof}
For each row-finite graph $E$, there is a unique monoid
homomorphism $\gamma_E\colon M_E\to V(L(E))$ such that $\gamma _E
(a_v)=[p_v]$. Clearly this defines a natural transformation from
the functor $M$ to the functor $V\circ L$; that is, if $f\colon
E\to F $ is a complete graph homomorphism, then the following
diagram commutes
\begin{equation*}\begin{CD}\notag
M_E @>{\gamma _E}>> V(L(E))\\
@V{M(f)}VV @VV{V(L(f))}V \\
M_F  @>{\gamma _F}>> V(L(F))
\end{CD} \notag \end{equation*}
We need to show that $\gamma_E$ is a monoid isomorphism for every
row-finite graph $E$. By using Lemma \ref{mondirectlimit} and
Lemma \ref{algdirectlimit}, we see that it is enough to show that
$\gamma _E$ is an isomorphism for a finite graph $E$.

Let $E$ be a finite graph and assume that $\{v_1,\dots
,v_m\}\subseteq E^0$ is the set of vertices which emit edges. We
start with an algebra
$$A_0=\prod _{v\in E^0}K .$$
In $A_0$ we have a family $\{ p_v:v\in E^0\}$ of orthogonal
idempotents such that $\sum_{v\in E^0} p_v=1$. Let us consider the
two finitely generated projective left $A_0$-modules
$P=A_0p_{v_1}$ and $Q=\oplus _{\{e\in E^1\mid
s(e)=v_1\}}A_0p_{r(e)}$. There exists an algebra $A_1:= A_0\langle
i,i^{-1}\colon \overline{P}\cong \overline{Q}\rangle $ with a
universal isomorphism $i\colon \overline{P}:=A_1\otimes _{A_0}P\to
\overline{Q}:=A_1\otimes _{A_0}Q$, see \cite[page 38]{Bergman}.
Note that this algebra is precisely the algebra $L(X_1)$, where
$X_1$ is the graph having $X_1^0=E^0$, and where $v_1$ emits the
same edges as it does in $E$, but all other vertices do not emit
any edge. Namely the row $(x_e: s(e)=v_1)$ implements an
isomomorphism $\overline{P}=A_1p_{v_1}\to \overline{Q}=\oplus
_{\{e\in E^1\mid s(e)=v_1\}}A_1p_{r(e)}$ with inverse given by the
column $(y_e:s(e)=v_1)^T$, which is clearly universal. By
\cite[Theorem 5.2]{Bergman}, the monoid $V(A_1)$ is obtained from
$V(A_0)$ by adjoining the relation $[P]=[Q]$. In our case we have
that $V(A_0)$ is the free abelian group on generators $\{a_v\mid
v\in E^0\}$, where $a_v=[p_v]$, and so $V(A_1)$ is given by
generators $\{a_v\mid v\in E^0\}$ and a single relation
$$a_{v_1}= \sum _{\{e\in E^1\mid s(e)=v_1\}}a_{r(e)}.$$

Now we proceed inductively. For $k\ge 1$, let $A_k$ be the graph
algebra $A_k=L(X_k)$, where $X_k$ is the graph with the same
vertices as $E$, but where only the first $k$ vertices $v_1,\dots
,v_k$ emit edges, and these vertices emit the same edges as they
do in $E$. Then we assume by induction that $V(A_k)$ is the
abelian group given by generators $\{a_v\mid v\in E^0\}$ and
relations
$$a_{v_i}= \sum _{\{e\in E^1\mid s(e)=v_i\}}a_{r(e)},$$
for $i=1,\dots ,k$. Let $A_{k+1}$ be the similar graph,
corresponding to vertices $v_1,\dots,v_k,v_{k+1}$. Then we have
$A_{k+1}=A_k\langle i,i^{-1}\colon \overline{P}\cong
\overline{Q}\rangle$ for $P=A_kp_{v_{k+1}}$ and $Q=\oplus_{\{e\in
E^1\mid s(e)=v_{k+1} \}}A_kp_{r(e)}$, and so we can apply again
Bergman's Theorem \cite[Theorem 5.2]{Bergman} to deduce that
$V(A_{k+1})$ is the monoid with the same generators as before and
the relations corresponding to $v_1,\dots,v_k,v_{k+1}$. It also
follows from \cite[Theorem 5.2]{Bergman} that the global dimension
of $L(E)$ is $\le 1$. This concludes the proof.
\end{proof}

\begin{example}
\label{example} {\rm Consider the following graph $E$:
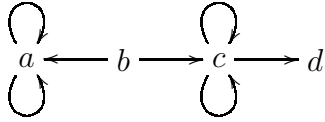
\[
{
\def\labelstyle{\displaystyle}
\xymatrix{ a\uloopr{}\dloopr{} & b \ar[l]\ar[r] & c
\uloopr{}\dloopr{} \ar[r] & d }}
\]
Then $M_E$ is the monoid generated by $a,b,c,d$ with defining
relations $a=2a$, $b=a+c$, $c=2c+d$. The Grothendieck group of
$M_E$ is infinite cyclic generated by the class of $c$. It follows
that $K_0(L_K(E))$ is infinite cyclic generated by $[p_c]$, and
$K_0(L_K(E))=K_0(L_K(E))^+$.}
\end{example}

\section{Refinement}

In this section we begin our formal study of the monoid $M_E$
associated with a row-finite graph $E$, and we show that $M_E$ is
a refinement monoid. The main tool is a careful description of the
congruence on the free abelian monoid given by the defining
relations of $M_E$.

Let $F$ be the free abelian monoid on the set $E^0$. The nonzero
elements of $F$ can be written in a unique form up to permutation
as $\sum _{i=1}^n x_i$, where $x_i\in E^0$. Now we will give a
description of the congruence on $F$ generated by the relations
(\ref{(M)}) on $F$. It will be convenient to introduce the
following notation. For $x\in E^0$, write
$${\bf r}(x):=\sum _{\{e\in E^1\mid s(e)=x\}} r(e)\in F .$$
With this new notation relations (\ref{(M)}) become $x={\bf r}(x)$
for every $x\in E^0$ that emits edges.

\begin{definition}
\label{binary} {\rm Define a binary relation $\rightarrow_1$ on
$F\setminus \{0\}$ as follows. Let $\sum _{i=1}^n x_i$ be an
element in $F$ as above and let $j\in \{1,\dots ,n\}$ be an index
such that $x_j$ emits edges. Then $\sum _{i=1}^n x_i\rightarrow_1
\sum _{i\ne j}x_i+{\bf r}(x_j)$. Let $\rightarrow $ be the
transitive and reflexive closure of $\rightarrow _1$ on
$F\setminus \{0\}$, that is, $\alpha\rightarrow \beta$ if and only
if there is a finite string $\alpha =\alpha _0\rightarrow _1
\alpha _1\rightarrow _1 \cdots \rightarrow _1 \alpha _t=\beta.$
Let $\sim$ be the congruence on $F$ generated by the relation
$\rightarrow_1 $ (or, equivalently, by the relation $\rightarrow
$). Namely $\alpha\sim \alpha$ for all $\alpha \in F$ and, for
$\alpha,\beta \ne 0$, we have $\alpha\sim \beta$ if and only if
there is a finite string $\alpha =\alpha _0,\alpha_1, \dots
,\alpha _n=\beta$, such that, for each $i=0,\dots ,n-1$, either
$\alpha _i\rightarrow _1 \alpha _{i+1}$ or
$\alpha_{i+1}\rightarrow_1 \alpha _i$. The number $n$ above will
be called the {\it length} of the string.}\qed
\end{definition}

It is clear that $\sim $ is the congruence on $F$ generated by
relations (\ref{(M)}), and so $M_E=F/{\sim}$.

The {\em support} of an element $\gamma$ in $F$, denoted
$\mbox{supp}(\gamma)\subseteq E^0$, is the set of basis elements
appearing in the canonical expression of $\gamma$.

\begin{lemma}
\label{division} Let $\rightarrow $ be the binary relation on $F$
defined above. Assume that $\alpha =\alpha _1+\alpha _2$ and
$\alpha \rightarrow \beta $. Then $\beta $ can be written as
$\beta =\beta _1+\beta _2$ with $\alpha _1\rightarrow \beta _1$
and $\alpha_2\rightarrow \beta _2$.
\end{lemma}

\begin{proof}
By induction, it is enough to show the result in the case where
$\alpha \rightarrow _1\beta$. If $\alpha \rightarrow _1 \beta$,
then there is an element $x$ in the support of $\alpha$ such that
$\beta = (\alpha -x)+{\bf r}(x)$. The element $x$ belongs either
to the support of $\alpha _1$ or to the support of $\alpha _2$.
Assume, for instance, that the element $x$ belongs to the support
of $\alpha _1$. Then we set $\beta _1=(\alpha _1-x)+{\bf r}(x)$
and $\beta _2=\alpha _2$.
\end{proof}

Note that the elements $ \beta _1$ and $\beta_2$ in Lemma
\ref{division} are not uniquely determined by $\alpha _1$ and
$\alpha _2$ in general, because the element $x\in E^0$ considered
in the proof could belong to {\it both} the support of $\alpha _1$
and the support of $\alpha _2$.

\smallskip

The following lemma gives the important ``confluence" property of
the congruence $\sim$ on the free abelian monoid $F$.

\begin{lemma}
\label{confluence} Let $\alpha$ and $\beta$ be nonzero elements in
$F$. Then $\alpha \sim \beta$ if and only if there is $\gamma\in
F$ such that $\alpha \rightarrow \gamma$ and $\beta \rightarrow
\gamma$.
\end{lemma}
\begin{proof}
Assume that $\alpha \sim \beta$. Then there exists a finite string
$\alpha =\alpha _0,\alpha_1, \dots ,\alpha _n=\beta$, such that,
for each $i=0,\dots ,n-1$, either $\alpha _i\rightarrow _1 \alpha
_{i+1}$ or $\alpha_{i+1}\rightarrow_1 \alpha _i$. We proceed by
induction on $n$. If $n=0$, then $\alpha=\beta$ and there is
nothing to prove. Assume the result is true for strings of length
$n-1$, and let $\alpha =\alpha_0,\alpha_1, \dots ,\alpha_n=\beta$
be a string of length $n$. By induction hypothesis, there is
$\lambda\in F$ such that $\alpha\rightarrow \lambda$ and
$\alpha_{n-1}\rightarrow \lambda$. Now there are two cases to
consider. If $\beta\rightarrow_1 \alpha_{n-1}$, then $\beta
\rightarrow \lambda$ and we are done. Assume that
$\alpha_{n-1}\rightarrow _1 \beta$. By definition of $\rightarrow
_1$, there is a basis element $x\in E^0$ in the support of $\alpha
_{n-1}$ such that $\alpha _{n-1}=x+\alpha_{n-1}'$ and $\beta =
{\bf r}(x)+\alpha _{n-1}'$. By Lemma \ref{division}, we have
$\lambda = \lambda (x)+\lambda '$, where $x\rightarrow \lambda
(x)$ and $\alpha_{n-1}'\rightarrow \lambda '$. If the length of
the string from $x$ to $\lambda (x)$ is positive, then we have
${\bf r}(x)\rightarrow \lambda (x)$ and so $\beta={\bf r}(x)
+\alpha _{n-1}'\rightarrow \lambda (x)+\lambda '=\lambda$. In case
that $x=\lambda (x)$, then set $\gamma= {\bf r}(x) +\lambda '$.
Then we have $\lambda \rightarrow _1 \gamma $ and so $\alpha
\rightarrow \gamma $, and also $\beta ={\bf r}(x) +\alpha
_{n-1}'\rightarrow {\bf r}(x)+\lambda '=\gamma $. This concludes
the proof.
\end{proof}

We are now ready to show the refinement property of $M_E$.

\begin{proposition}
\label{refinement} The monoid $M_E$ associated with any row-finite
graph $E$ is a refinement monoid.
\end{proposition}
\begin{proof}
Let $\alpha=\alpha _1+\alpha_2\sim\beta =\beta _1+\beta _2$, with
$\alpha_1,\alpha_2,\beta _1,\beta _2\in F$. By Lemma
\ref{confluence}, there is $\gamma\in F$ such that $\alpha
\rightarrow \gamma$ and $\beta \rightarrow \gamma$. By Lemma
\ref{division}, we can write $\gamma =\alpha_1 '+\alpha _2'=\beta
_1'+\beta _2'$, with $\alpha_i\rightarrow \alpha _i'$ and $\beta
_i\rightarrow \beta _i'$ for $i=1,2$. Since $F$ is a free abelian
monoid, $F$ has the refinement property and so there are
decompositions $\alpha _i'=\gamma _{i1}+\gamma _{i2}$ for $i=1,2$
such that $\beta _j'=\gamma _{1j}+\gamma _{2j}$ for $j=1,2$. The
result follows.
\end{proof}

\section{Ideal lattice}

Let $E=(E^0,E^1)$ be a row-finite directed graph. In this section,
we will establish the connection between saturated hereditary
subsets of $E^0$, order-ideals of the associated monoid $M_E$, and
graded ideals of the graph algebra $L_K(E)$.

We start by recalling some basic concepts in graph theory, that
will be needed later.

Let $E=(E^0,E^1)$ be a directed graph. For $n\ge 2$, we define
$$E^n:=\left\{\alpha=(\alpha_1,\dots,\alpha_n)\mid \, \alpha_i \in E^1
\mbox{ and } r(\alpha_i)=s(\alpha_{i+1}) \mbox{ for } 1\leq i\leq
n-1\right\},
$$
and $E^*=\bigcup_{n\ge 0}E^n$.

We define a relation $\geq$ on $E^0$ by setting $v\geq w$ if there
is a path $\mu \in E^*$ with $s(\mu)=v$ and $r(\mu)=w$. A subset
$H$ of $E^0$ is called {\em hereditary} if $v \geq w$ and $v\in H$
imply $w\in H$. A hereditary set $H$ is {\em saturated} if every
vertex which feeds into $H$ and only into $H$ is again in $H$;
that is, if $s^{-1}(v)\neq \emptyset$ and $r(s^{-1}(v))\subseteq
H$, then $v\in H$.

\begin{definition}\label{treeofv}\mbox{ }{\rm
Let $v \in E^0$. We define the tree of $v$, to be the subset of
$E^0$
$$T(v)=\left\{ w\in E^0\mid\, \exists \, \alpha \in E^* \mbox{ with }
s(\alpha)=v \mbox{ and }r(\alpha)=w \right\}=\{w\in E^0\mid v\ge
w\}.$$} \end{definition} Clearly, the tree of $v$ is the smallest
hereditary subset of $E^0$ containing $v$.

We denote by $\mathcal{H}$ the set of saturated hereditary subsets
of the graph $E$.

Since the intersection of saturated sets is saturated, there is a
smallest saturated subset $\overline S $ containing any given
subset $S$ of $E^0$. We will call $\overline S$ the saturation of
$S$. The saturation $\overline H$ of a hereditary set $H$ is again
hereditary. Indeed, $\overline H=\bigcup _{n=0}^{\infty} \Lambda
_n(H)$ is an increasing union of hereditary subsets $\Lambda
_n(H)$, for $n\ge 0$, which are defined inductively as follows:
\begin{enumerate}
\item $\Lambda _0(H)=H$. \item $\Lambda _n(H)=\{y\in E^0\mid
s^{-1}(y)\ne \emptyset \text{ and }r(s^{-1}(y))\subseteq \Lambda
_{n-1}(H)\} \cup \Lambda _{n-1}(H)$, for $n\ge 1$.
\end{enumerate}
In particular this applies to the hereditary subsets of the form
$T(v)$, where $v\in E^0$: The saturated hereditary subset of $E$
generated by $v$ is $\overline{T}(v)=\bigcup
_{n=0}^{\infty}\Lambda _n(T(v))$.

An {\em order-ideal} of a monoid $M$ is a submonoid $I$ of $M$
such that $x+y=z$ in $M$ and $z\in I$ imply that both $x,y$ belong
to $I$. An order-ideal can also be described as a submonoid $I$ of
$M$, which is hereditary with respect to the canonical pre-order
$\le $ on $M$: $x\le y$ and $y\in I$ imply $x\in I$. Recall that
the pre-order $\le $ on $M$ is defined by setting $x\le y$ if and
only if there exists $z\in M$ such that $y=x+z$.

The set $\mathcal{L}(M)$ of order-ideals of $M$ forms a (complete)
lattice $\left(\mathcal L(M), \subseteq, \overline{\sum}, \cap
\right)$. Here, for a family of order-ideals $\{I_i\}$, we denote
by $\overline{\sum }I_i$ the set of elements $x\in M$ such that
$x\le y$, for some $y$ belonging to the algebraic sum $\sum I_i$
of the order-ideals $I_i$. Note that $\sum I_i=\overline{\sum}I_i$
whenever $M$ is a refinement monoid.

Let $F_E$ be the free abelian monoid on $E^0$, and recall that
$M_E=F_E/{\sim}$. For $\gamma \in F_E$ we will denote by
$[\gamma]$ its class in $M_E$. Note that any order-ideal $I$ of
$M_E$ is generated {\em as a monoid} by the set $\{[v]\mid v\in
E^0\}\cap I$.

The set $\mathcal H$ of saturated hereditary subsets of $E^0$ is
also a complete lattice $\left(\mathcal{H}, \subseteq,
\overline{\cup}, \cap\right)$.

\begin{proposition}
\label{lattmonoid} Let $E$ be a row-finite graph. Then, there are
order-preserving mutually inverse maps
$$ \varphi\colon \mathcal{H} \longrightarrow
\mathcal L (M_E); \qquad \psi\colon\mathcal{L}(M_E)
\longrightarrow \mathcal{H},$$ where $\varphi (H)$ is the
order-ideal of $M_E$ generated by $\{[v]\mid v\in H\}$, for $H\in
\mathcal H$, and $\psi (I)$ is the set of elements $v$ in $E^0$
such that $[v]\in I$, for $I\in \mathcal L(M_E)$.
\end{proposition}

\begin{proof}
The maps $\varphi$ and $\psi$ are obviously order-preserving. It
will be enough to show the following facts:
\begin{enumerate}
\item For $I\in \mathcal L (M_E)$, the set $\psi (I)$ is a
hereditary and saturated subset of $E^0$. \item If $H\in \mathcal
H$ then $[v]\in \varphi(H)$ if and only if $v\in H$.
\end{enumerate}
For, if (1) and (2) hold true, then $\psi$ is well-defined by (1),
and  $\psi (\varphi (H))=H$ for $H\in \mathcal H$, by (2). On the
other hand, if $I$ is an order-ideal of $M_E$, then obviously
$\varphi (\psi (I))\subseteq I$, and since $I$ is generated as a
monoid by $\{[v]\mid v\in E^0\}\cap I=[\psi (I)]$, it follows that
$I\subseteq \varphi (\psi (I))$.

{\em Proof of (1):} Let $I$ be an order-ideal of $M_E$, and set
$H:=\psi (I)=\{v\in E^0\mid [v]\in I\}$. To see that $H$ is
hereditary, we have to prove that, whenever we have a path
$(e_1,e_2,\dots ,e_n)$ in $E$ with $s(e_1)=v$ and $r(e_n)=w$ and
$v\in H$, then $w\in H$. If we consider the corresponding path
$v\to _1 \gamma_1\to _1 \gamma _2\to _1\cdots \to _1 \gamma _n$ in
$F_E$, we see that $w$ belongs to the support of $\gamma _n$, so
that $w\le \gamma _n$ in $F_E$. This implies that $[w]\le [\gamma
_n]=[v]$, and so $[w]\in I$ because $I$ is hereditary.

To show saturation, take $v$ in $E^0$ such that $r(e)\in H$ for
every $e\in E^1$ such that $s(e)=v$. We then have
$\mbox{supp}({\bf r}(v))\subseteq H$, so that $[{\bf r}(v)]\in I$
because $I$ is a submonoid of $M_E$. But $[v]=[{\bf r}(v)]$, so
that $[v]\in I$ and $v\in H$.

{\em Proof of (2):} Let $H$ be a saturated hereditary subset of
$E^0$, and let $I:=\varphi (H)$ be the order-ideal of $M_E$
generated by $\{[v]\mid v\in H\}$. Clearly $[v]\in I$ if $v\in H$.
Conversely, suppose that $[v]\in I$. Then $[v]\le [\gamma]$, where
$\gamma \in F_E $ satisfies $\mbox{supp}(\gamma)\subseteq H $.
Thus we can write $[\gamma] =[v]+[\delta]$ for some $\delta\in
F_E$. By Lemma \ref{confluence}, there is $\beta \in F_E$ such
that $\gamma\to \beta$ and $v+\delta \to \beta$. Since $H$ is
hereditary and $\mbox{supp}(\gamma)\subseteq H$, we get
$\mbox{supp}(\beta )\subseteq H$. By Lemma \ref{division}, we have
$\beta =\beta _1+\beta _2$, where $v\to \beta _1$ and $\delta \to
\beta _2$. Observe that $\mbox{supp}(\beta _1)\subseteq
\mbox{supp}(\beta)\subseteq H$. Using that $H$ is saturated, it is
a simple matter to check that, if $\alpha \to _1 \alpha '$ and
$\mbox{supp}(\alpha ')\subseteq H$, then $\mbox{supp}(\alpha
)\subseteq H$. Using this and induction, we obtain that $v\in H$,
as desired.
\end{proof}

We next consider ideals in the algebra $L_K(E)$ associated with
the graph $E$. For a general unital ring $R$, the lattice of
order-ideals of $V(R)$ is isomorphic with the lattice of trace
ideals of $R$; see \cite{AF} and \cite{FHK}. It is straightforward
to see that this lattice isomorphism also holds when $R$ is a ring
with local units. In particular, the lattice of order-ideals of
$V(L_K(E))$ is isomorphic with the lattice of trace ideals of
$L_K(E)$. Being $V(L_K(E))\cong M_E $ a refinement monoid
(Proposition \ref{refinement}), we see that the trace ideals of
$L_K(E)$ are exactly the ideals generated by idempotents of
$L_K(E)$. In general not all the ideals in $L_K(E)$ will be
generated by idempotents. For instance, if $E$ is a single loop,
then $L_K(E)=K[x,x^{-1}]$ and the ideal generated by $1-x$ only
contains the idempotent $0$. However, it is possible to describe
the ideals generated by idempotents by using the canonical grading
of $L_K(E)$. Let us recall that $L_K(E)$ is generated by sets
$\{p_v\mid v\in E^0\}$ and $\{x_e,y_e\mid e\in E^1\}$, which
satisfy the following relations:

(1) $p_vp_{v'}=\delta _{v,v'}p_v$ for all $v,v'\in E^0$.

(2) $p_{s(e)}x_e=x_ep_{r(e)}=x_e$ for all $e\in E^1$.

(3) $p_{r(e)}y_e=y_ep_{s(e)}=y_e$ for all $e\in E^1$.

(4) $y_ex_{e'}=\delta _{e,e'}p_{r(e)}$ for all $e,e'\in E^1$.

(5) $p_v=\sum _{\{ e\in E^1\mid s(e)=v \}}x_ey_e$ for every $v\in
E^0$ that emits edges.

If we declare that the degree of $x_e$ is $1$ and the degree of
$y_e$ is $-1$ for all $e\in E^1$, and that the degree of each
$p_v$ is $0$ for $v\in E^0$, then we obtain a well-defined degree
on the algebra $L(E)=L_K(E)$, because all relations (1)--(5) are
homogeneous. Thus $L(E)$ is a $\Z$-graded algebra:
$$L(E)=\bigoplus _{n\in \Z} L(E)_n ;
\qquad  L(E)_nL(E)_m\subseteq L(E)_{n+m}, ,\,\ \text{ for all }
n,m\in \Z. $$ For a subset $X$ of a $\Z$-graded ring $R=\oplus
_{n\in \Z} R_n$, set $X_n=X\cap R_n$. An ideal $I$ of $R$ is said
to be a {\em graded ideal} in case $I=\bigoplus _{n\in \Z}I_n$.
Let us denote the lattice of graded ideals of a $\Z$-graded ring
$R$ by $\mathcal L _{\text{gr}}(R)$.

Recall that $v\in E^0$ is called a {\em sink} in case $v$ does not
emit any edge.

\begin{theorem}
\label{lattall} Let $E$ be a row-finite graph. Then there are
order-isomorphisms $$\mathcal H\cong \mathcal L (M_E)\cong
\mathcal L_{\text{gr}}(L_K(E)),$$ where $\mathcal H $ is the
lattice of hereditary and saturated subsets of $E^0$, $\mathcal L
(M_E)$ is the lattice of order-ideals of the monoid $M_E$, and
$\mathcal L_{\text{gr}}(L_K(E))$ is the lattice of graded ideals
of the graph algebra $ L_K(E)$.
\end{theorem}

\begin{proof}
We have obtained an order-isomorphism $\mathcal H\cong \mathcal L
(M_E)$ in Proposition \ref{lattmonoid}. As we observed earlier
there is an order-isomorphism $\mathcal L(M_E)=\mathcal L
(V(L(E)))\cong \mathcal L_{\text{idem}}(L(E))$, where $\mathcal
L_{\text{idem}}(L(E))$ is the lattice of ideals in $L(E)$
generated by idempotents. The isomorphism is given by the rule
$I\mapsto \widetilde{I}$, for every order-ideal $I$ of $M_E$,
where $\widetilde{I}$ is the ideal generated by all the
idempotents $e\in L(E)$ such that $V(e)\in I$. (Here $V(e)$
denotes the class of $e$ in $V(L(E))=M_E$.) Given any order-ideal
$I$ of $M_E$, it is generated as monoid by the elements
$V(p_v)(=[v]=a_v )$ such that $V(p_v)\in I$, so that
$\widetilde{I}$ is generated as an ideal by the idempotents $p_v$
such that $p_v\in \widetilde{I}$. In particular we see that every
ideal of $L(E)$ generated by idempotents is a graded ideal.

It only remains to check that every graded ideal of $L(E)$ is
generated by idempotents. For this, it will be convenient to
recall the definition of the path algebra $P(E)$ associated with
$E$. The algebra $P(E)$ is the algebra generated by a set
$\{p_v\mid v\in E^0\}$ of pairwise orthogonal idempotents,
together with a set of variables $\{x_e\mid e\in E^1\}$, which
satisfy relation (2). A $K$-basis for $P(E)$ is given by the set
of ``paths" $\gamma =x_{e_1}x_{e_2}\cdots x_{e_r}$, such that
$r(e_i)=s(e_{i+1})$ for $i=1,\dots ,r-1$. We put
$s(\gamma)=s(e_1)$ and $r(\gamma)=r(e_r)$, and the {\em length}
$|\gamma |$ of $\gamma$ is defined to be $r$. It is easy to see
that $P(E)$ is indeed a subalgebra of $L(E)$. The algebra $P(E)^*$
is, by definition, the subalgebra of $L(E)$ generated by
$\{p_v\mid v\in E^0\}$ and $\{y_e\mid e\in E^1\}$. Of course we
can define an involution on $L(E)$ sending $x_e$ to $y_e$, so that
all $p_v$ are projections: $p_v=p_v^2=p_v^*$, and acting on $K$ by
any prescribed involution on $K$. Note that for $\gamma
=x_{e_1}x_{e_2}\cdots x_{e_r}$, we have $\gamma ^*=y_{e_r}\cdots
y_{e_2}y_{e_1}$. Now elements in $L(E)$ can be described as linear
combinations of elements of the form $\gamma \nu ^*$, where
$\gamma$ and $\nu $ are paths on $E$ with $r(\gamma)=r(\nu)$. It
is clear that, for $n>0$, we have $L(E)_n=\oplus _{|\gamma |=n}
\gamma L(E)_0$, and similarly, $L(E)_{-n}=\oplus
_{|\gamma|=n}L(E)_0 \gamma ^*$.

Given a graded ideal $J$ of $L(E)$, take any element $a\in J_n$,
where $n>0$. Then $a=\sum _{|\gamma |=n}\gamma a_{\gamma}$, for
some $a_{\gamma}\in L(E)_0$. For a fixed path $\nu$ of length $n$,
we have $\nu ^*a=a_{\nu}$, so that $a_{\nu}\in J_0$. We conclude
that $J_n=L(E)_nJ_0$, and similarly $J_{-n}=J_0L(E)_{-n}$. Since
$J$ is a graded ideal, we infer that $J$ is generated as ideal by
$J_0$, which is an ideal of $L(E)_0$.

To conclude the proof, we only have to check that every ideal of
$L(E)_0$ is generated by idempotents. Indeed we will prove that
$L(E)_0$ is a von Neumann regular ring, more precisely $L(E)_0$ is
an ultramatricial $K$-algebra, i.e. a direct limit of matricial
algebras over $K$ \cite{vnrr}, though not all the connecting
homomorphisms are unital. (A matricial $K$-algebra is a finite
direct product of full matrix algebras over $K$.)

By Lemma \ref{algdirectlimit} we have $L(E)=\varinjlim _{i\in
I}L(X_i)$ for a directed family $\{X_i\mid i\in I\}$ of finite
graphs. Then $L(E)_0=\varinjlim_{i\in I} L(X_i)_0$, and so we can
assume that $E$ is a finite graph.

Now for a finite graph $E$, all the transition maps are unital.
They can be built in the following fashion. For each $v$ in $E^0$,
and each $n\in \Z^+$, let us denote by $P(n,v)$ the set of paths
$\gamma= x_{e_1}\cdots x_{e_n}\in P(E)$ such that $|\gamma |=n$
and $r(\gamma)=v$.  The set of sinks will be denoted by $S(E)$.
Now the algebra $L(E)_0$ admits a natural filtration by algebras
$L_{0,n}$, for $n\in \Z^+$. Namely $L_{0,n}$ is the set of linear
combinations of elements of the form $\gamma \nu^*$, where
$\gamma$ and $\nu$ are paths with $r(\gamma)=r(\nu)$ and
$|\gamma|=|\nu|\le n$. The algebra $L_{0,0}$ is isomorphic to
$\prod _{v\in E^0}K$. In general the algebra $L_{0,n}$ is
isomorphic to
$$\big[ \prod _{i=0}^{n-1}\big( \prod _{v\in S(E)}M_{|P(i,v)|}(K)\big)\big]
\times \big[ \prod _{v\in E^0}M_{|P(n,v)|}(K) \big] . $$ The
transition homomorphism $L_{0,n}\to L_{0,n+1}$ is the identity on
the factors $\prod _{v\in S(E)}M_{|P(i,v)|}(K)$, for $0\le i\le
n-1$, and also on the factor $\prod _{v\in S(E)}M_{|P(n,v)|}(K)$
of the last term of the displayed formula. The transition
homomorphism
$$\prod_{v\in E^0\setminus S(E)}M_{|P(n,v)|}(K)\to \prod_ {v\in
E^0}M_{|P(n+1,v)|}(K)$$ is a block diagonal map induced by the
following identification in $L(E)_0$: A matrix unit in a factor
$M_{|P(n,v)|}(K)$, where $v\in E^0\setminus S(E)$, is a monomial
of the form $\gamma\nu ^*$, where $\gamma$ and $\nu$ are paths of
length $n$ with $r(\gamma)=r(\nu)=v$. Since $v$ is not a sink, we
can enlarge the paths $\gamma$ and $\nu$ using the edges that $v$
emits, obtaining paths of length $n+1$, and relation (5) in the
definition of $L(E)$ gives $\gamma \nu ^*=\sum _{\{e\in E^1\mid
s(e)=v\}}(\gamma x_e)(y_e\nu^*)$.

It follows that $L(E)_0$ is an ultramatricial $K$-algebra, and the
proof is complete.
\end{proof}

\section{Separativity}

In this section we prove that the monoid $M_E$ associated with a
row-finite graph $E=(E^0,E^1)$ is always a separative monoid.
Recall that this means that for elements $x,y,z\in M_E$, if
$x+z=y+z$ and $z\le nx$ and $z\le ny$ for some positive integer
$n$, then $x=y$.

The separativity of $M_E$ follows from results of Brookfield
\cite{Brook} on primely generated monoids; see also \cite[Chapter
6]{Wehrung}. Indeed the class of primely generated refinement
monoids satisfies many other nice cancellation properties. We will
highlight unperforation later, and refer the reader to
\cite{Brook} for further information.

\begin{definition}
\label{prime} {\rm Let $M$ be a monoid. An element $p\in M$ is
prime if for all $a_1,a_2\in M$, $p\le a_1+a_2$ implies $p\le a_1$
or $p\le a_2$. A monoid is primely generated if each of its
elements is a sum of primes.}
\end{definition}

\begin{proposition} \cite[Corollary 6.8]{Brook}
\label{fingen} Any finitely generated refinement monoid is primely
generated.
\end{proposition}

It follows from Proposition \ref{fingen} that, for a finite graph
$E$, the monoid $M_E$ is primely generated. Note that this is not
always the case for a general row-finite graph $E$. An example is
provided by the graph:

$$\xymatrix{p_0\ar[r] \ar[d]  & p_1 \ar[r]\ar[dl] & p_2 \ar[r]\ar[dll] & p_3 \ar[r]\ar[dlll] & \cdots \\
a}$$

The corresponding monoid $M$ has generators $a,p_0,p_1, \dots $
and relations given by $p_i=p_{i+1}+a$ for all $i\ge 0$. One can
easily see that the only prime element in $M$ is $a$, so that $M$
is not primely generated.

\begin{theorem}
\label{separativity} Let $E$ be a row-finite graph. Then the
monoid $M_E$ is separative.
\end{theorem}

\begin{proof}
By Lemma \ref{mondirectlimit}, we get that $M_E$ is the direct
limit of monoids $M_{X_i}$ corresponding to finite graphs $X_i$.
Therefore, in order to check separativity, we can assume that the
graph $E$ is finite.

Assume that $E$ is a finite graph. Then $M_E$ is generated by the
finite set $E^0$ of vertices of $E$, and thus $M_E$ is finitely
generated. By Proposition \ref{refinement}, $M_E$ is a refinement
monoid, so it follows from Proposition \ref{fingen} that $M_E$ is
a primely generated refinement monoid. By \cite[Theorem
4.5]{Brook}, the monoid $M_E$ is separative.
\end{proof}

As we remarked before, primely generated refinement monoids
satisfy many nice cancellation properties, as shown in
\cite{Brook}. Some of these properties are preserved in direct
limits, so they are automatically true for the graph monoids
corresponding to any row-finite graph. Especially important in
several applications is the property of unperforation. Let us say
that a monoid $M$ is {\it unperforated} in case, for all elements
$a,b\in M$ and all positive integers $n$, we have $na\le
nb\implies a\le b$. This implies that the Grothendieck group
$G(M)$ of $M$ is {\it unperforated}: for all $g\in G(M)$ and all
positive integers $n$, we have $ng\ge 0\implies g\ge 0$.

\begin{proposition}
\label{unperforation} Let $E$ be a row-finite graph. Then the
monoid $M_E$ is unperforated.
\end{proposition}

\begin{proof}
As in the proof of Theorem \ref{separativity}, we can reduce to
the case of a finite graph $E$. In this case, the result follows
from \cite[Corollary 5.11(5)]{Brook}.
\end{proof}

\begin{corollary}
\label{separativerings} Let $E$ be a row-finite graph. Then
$FP(L_K(E))$ satisfies the refinement property and $L_K(E)$ is a
separative ring. Moreover, the monoid $V(L_K(E))$ is an
unperforated monoid and $K_0(L_K(E))$ is an unperforated group.
\end{corollary}

\begin{proof}
By Theorem \ref{VLE}, we have $V(L_K(E))\cong M_E$. So the result
follows from Proposition \ref{refinement}, Theorem
\ref{separativity} and Proposition \ref{unperforation}.
\end{proof}

\smallskip

Another useful technique to deal with graph monoids of finite
graphs consists in considering composition series of order-ideals
in the monoid. These composition series correspond via Theorem
\ref{lattall} to composition series of graded ideals in $L_K(E)$,
and, using \cite[Theorem 4.1(b)]{BPRS}, they also correspond to
composition series of closed gauge-invariant ideals of the graph
$C^*$-algebra $C^*(E)$. This approach will be used in the proof of
Theorem \ref{C*algebras} in our next section. It also leads to a
different proof of the separativity of $M_E$ (Theorem
\ref{separativity}), that will be sketched in Remark
\ref{remarkseparativity}.

Given an order-ideal $S$ of a monoid $M$ we define a congruence
$\sim _S$ on $M$ by setting $a\sim _S b$ if and only if there
exist $e,f\in S$ such that $a+e=b+f$. Let $M/S$ be the factor
monoid obtained from the congruence $\sim _S$; see \cite{AGOP}.
For large classes of rings $R$, one has $V(R/I)\cong V(R)/V(I)$
for any ideal $I$ of $R$; see \cite[Proposition 1.4]{AGOP}.
\smallskip

We need a monoid version of \cite[Theorem 4.1(b)]{BPRS}.

\begin{lemma}
\label{factoring} Let $E$ be a row-finite graph. For a saturated
hereditary subset $H$ of $E^0$, consider the order-ideal
$S=\varphi (H)$ associated with $H$, as in Proposition
\ref{lattmonoid}. Let $G=(G^0,G^1)$ be the graph defined as
follows. Put $G^0= E^0\setminus H$ and $G^1=\{e\in E^1\mid r(e)\in
G^0\}$. Then there is a natural monoid isomorphism $M_E/S\cong
M_G$.
\end{lemma}

\begin{proof}
Note that $S$ is generated as a monoid by the elements $a_v$, with
$v\in H$. There is a unique monoid homomorphism $\pi \colon F_E\to
F_G$ sending $v$ to $0$ for $v\in H$ and $v$ to $v$ for $v\in
E\setminus H$, where $F_E$ (respectively $F_G$) is the free
abelian monoid on $E$ (respectively $G$). The map $\pi $ induces a
surjective monoid homomorphism $\ol{\pi}\colon M_E\to M_G $, and
it is clear that $\ol{\pi} $ factors through $M_E/S$, i.e, we have
$M_E\to M_E/S\to M_G$. If $\pi (\alpha)\sim \pi (\beta)$ in $F_G$
for $\alpha,\beta \in F_E$, then by Lemma \ref{confluence} there
is $\gamma\in F_G$ such that $\pi(\alpha)\to \gamma$ and $\pi
(\beta)\to \gamma$. This means that there is a string
$\pi(\alpha)=\gamma_0\to _1\gamma_1\to _1\cdots \to
_1\gamma_r=\gamma$ in $F_G$, and similarly for $\pi(\beta)\to
\gamma$. Let us consider the same strings, but now in $F_E$. We
then get that $\alpha\to \gamma+\delta_1$, where $\delta _1$ is
supported on $H$, and similarly $\beta\to \gamma+\delta _2$, where
$\delta _2$ is supported on $H$. It follows that the following
identity holds in $M_E$:
$$[\alpha] +[\delta _2]=[\gamma]+[\delta _1]+[\delta _2]=[\beta]+[\delta _1],$$
with $[\delta_1],[\delta _2]\in S$. We conclude that the map
$M_E/S\to M_G$ is injective, and so it is a monoid isomorphism, as
desired.
\end{proof}

Let us call $M_E$ {\it simple} if $M_E$ has only the trivial
order-ideals. This corresponds by Proposition \ref{lattmonoid} to
the situation where the hereditary and saturated subset generated
by any vertex of $E$ is $E^0$. It is well-known that this happens
if and only if $E$ is cofinal. Let $E^{\le \infty}$ be the set of
infinite paths in $E$ together with the finite paths in $E$ whose
end point is a sink. Then $E$ is said to be {\em cofinal} in case
given a vertex $v$ in $E$ and a path $\gamma$ in $E^{\le \infty}$,
there is a vertex $w$ in the path $\gamma$ such that $v\ge w$.

A finite path $\alpha$ of positive length is called a {\em loop}
if $s(\alpha)=r(\alpha)=v$. A loop $\alpha=(e_1,e_2, \dots ,e_n)$
is {\em simple} if all the vertices $s(e_i)$, $1\le i\le n$, are
distinct. For a subgraph $G$ of $E$, an {\em exit} of $G$ is an
edge $e$ in $E$ with $s(e)\in G^0$ and $e\notin G^1$.

\begin{remark}
\label{remarkseparativity} {\rm We are now ready to sketch a
different proof of the separativity of $M_E$ (Theorem
\ref{separativity}), using the theory of order-ideals.

As in the proof of Theorem \ref{separativity}, we can assume that
$E$ is a finite graph. In this case it is obvious that $E^0$ has a
finite number of saturated hereditary subsets, so $M_E$ has a
finite number of order-ideals. Take a finite chain $0=S_0\le
S_1\le \cdots \le S_n=M_E$ such that each $S_i$ is an order-ideal
of $M_E$, and all the quotients $S_i/S_{i-1}$ are simple. By
Proposition \ref{lattmonoid}, we have $S_i\cong M_{H_i}$, for some
finite graph $H_i$, and by Lemma \ref{factoring}, we have
$S_i/S_{i-1}\cong M_{G_i}$ for some cofinal finite graph $G_i$. By
Proposition \ref{refinement}, $S_i$ is a refinement monoid for all
$i$, so the Extension Theorem for refinement monoids
(\cite[Theorem 4.5]{AGOP}) tells us that $S_i$ is separative if
and only if so are $S_{i-1}$ and $S_i/S_{i-1}$. It follows by
induction that it is enough to show the case where $E$ is a
cofinal finite graph.

\smallskip

Let $E$ be a cofinal finite graph. We distinguish three cases.
First, suppose that $E$ does not have loops. Then there is a sink
$v$, and by cofinality for every vertex $w$ of $E$ there is a path
from $w$ to $v$. It follows that $M_E$ is a free abelian monoid of
rank one (i.e. isomorphic to $\Z^+$), generated by $a_v$. In
particular $M_E$ is a separative monoid. Secondly, assume that $E$
has a simple loop without exit, and let $v$ be any vertex in this
simple loop. By using the cofinality condition, it is easy to see
that there are no other simple loops in $E$, and that every vertex
in $E$ connects to $v$. It follows again that $M_E$ is a free
abelian monoid of rank one, generated by $a_v$.

Finally we consider the case where every simple loop has an exit.
By cofinality, every vertex connects to every loop. Using this and
the property that every loop has an exit, it is quite easy to show
that for every nonzero element $x$ in $M_E$ there is a nonzero
element $y$ in $M_E$ such that $x=x+y$. It follows that
$M_E\setminus \{0\}$ is a group; see for example \cite[Proposition
2.4]{AGOP}. In particular $M_E$ is a separative monoid.
\qed}\end{remark}

\begin{example}
\label{examplebis} {\rm We consider again the graph $E$ described
in Example \ref{example}. A composition series of order-ideals for
$M_E$ is obtained from the graph monoids corresponding to the
following chain of saturated hereditary subsets of $E$:
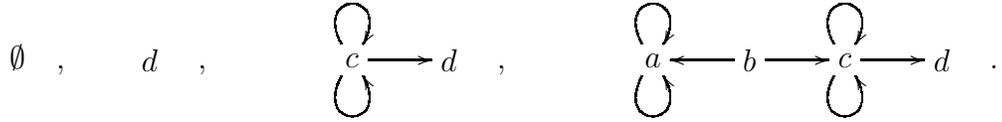
$$\emptyset \quad ,\qquad
 \xymatrix{
  d
}\quad , \qquad \qquad \xymatrix{
 c \uloopr{}\dloopr{}
\ar[r] & d } \quad ,\qquad \qquad \xymatrix{ a\uloopr{}\dloopr{} &
b \ar[l]\ar[r] & c \uloopr{}\dloopr{} \ar[r] & d }\quad .
$$
By Lemma \ref{factoring}, the corresponding simple quotient
monoids are the graph monoids corresponding to the following
graphs:
$$
 
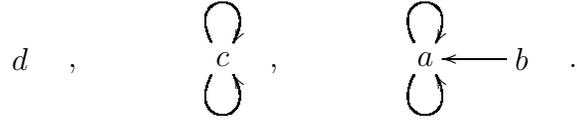
\xymatrix{
  d
}\quad , \qquad \qquad \xymatrix{
 c \uloopr{}\dloopr{}
} \quad ,\qquad \qquad \xymatrix{ a\uloopr{}\dloopr{} & b \ar[l]
}\quad .
$$
}
\end{example}

\section{The monoid associated with a graph $C^*$-algebra}

In this section, we will assume that $L(E)=L_{\mathbb C}(E)$ is
the graph algebra of the graph $E$ over the field $\C$ of complex
numbers, endowed with its natural structure of complex
$*$-algebra, so that $x_e^*=y_e$ for all $e\in E^1$, $p_v^*=p_v$
for all $v\in E^0$, and $(\xi a)^*=\overline{\xi}a^*$ for $\xi\in
\mathbb C$ and $a\in L(E)$. There is a natural inclusion of
complex $*$-algebras $\psi \colon L(E)\to C^*(E)$, where $C^*(E)$
denotes the graph $C^*$-algebra associated with $E$.

\begin{theorem}
\label{C*algebras} Let $E$ be a row-finite graph, and let
$L(E)=L_{\mathbb C}(E)$ be the graph algebra over the complex
numbers. Then the natural inclusion $\psi\colon L(E)\to C^*(E)$
induces a monoid isomorphism $V(\psi)\colon V(L(E))\to V(C^*(E))$.
In particular the monoid $V(C^*(E))$ is naturally isomorphic with
the monoid $M_E$.
\end{theorem}

\begin{proof}
The algebra homomorphism $\psi\colon L(E)\to C^*(E)$ induces the
following commutative square:
\begin{equation*}\begin{CD}\notag
V(L(E)) @>V(\psi)>> V(C^*(E))\\
@V{\varphi _1}VV @VV{\varphi _2}V \\
K_0(L(E))  @>K_0(\psi )>> K_0(C^*(E))
\end{CD} \notag \end{equation*}
The map $K_0(\psi)$ is an isomorphism by Theorem \ref{VLE} and
\cite[Theorem 3.2]{RS}. Using Lemma \ref{algdirectlimit} and Lemma
\ref{C*algdirectlimit}, we see that it is enough to show that
$V(\psi)$ is an isomorphism for a finite graph $E$.

Assume that $E$ is a finite graph. We first show that the map
$V(\psi )\colon V(L(E))\to V(C^*(E))$ is injective. Suppose that
$P$ and $Q$ are idempotents in $M_{\infty}(L(E))$ such that $P\sim
Q$ in $C^*(E)$. By Theorem \ref{VLE}, we can assume that each of
$P$ and $Q$ are equivalent in $M_{\infty}(L(E))$ to direct sums of
``basic" projections, that is, projections of the form $p_v$, with
$v\in E^0$. Let $J$ be the closed ideal of $C^*(E)$ generated by
the entries of $P$. Since $P\sim Q$, the closed ideal generated by
the entries of $P$ agrees with the closed ideal generated by the
entries of $Q$ and indeed it agrees with the closed ideal
generated by the projections of the form $p_w$, where $w$ ranges
on the saturated hereditary subset $H$ of $E^0$ generated by
$\{v\in E^0\mid P=\oplus p_v\}$ (see \cite[Theorem 4.1]{BPRS}). It
follows from Theorem \ref{lattall} that $P$ and $Q$ generate the
same ideal $I_0$ in $L(E)$. There is a projection $e\in L(E)$,
which is the sum of the basic projections $p_w$, where $w$ ranges
in $H$, such that $I_0=L(E)eL(E)$ and $eL(E)e=L(H)$ is also a
graph algebra. Note that $P$ and $Q$ are full projections in
$L(H)$, and so $[1_H]\le m[P]$ and $[1_H]\le m[Q]$ for some $m\ge
1$. Now consider the map $\psi _H\colon L(H)\to C^*(H)$. Since
$V(\psi _H)([P])= V(\psi_H)([Q]) $ in $V(C^*(H))$ we get
$K_0(\psi_H)(\varphi _1([P]))=K_0(\psi_H)(\varphi_1 ([Q]))$, and
since $K_0(\psi _H)$ is an isomorphism we get $\varphi _1
([P])=\varphi_1([Q])$. This means that there is $k\ge 0$ such that
$[P]+k[1_H]= [Q]+ k[1_H]$. But since $V(L(E))$ is separative and
$[1_H]\le m[P]$ and $[1_H]\le m[Q]$, we get $[P]=[Q]$ in
$V(L(E))$.

\smallskip

Now we want to see that the map $V(\psi)\colon V(L(E))\to
V(C^*(E))$ is surjective. By \cite[Theorem 4.1]{BPRS}, there is a
natural isomorphism between the lattice of saturated hereditary
subsets of $E^0$ and the lattice of closed gauge-invariant ideals
of $C^*(E)$. Thus, since $E$ is finite, the number of closed
gauge-invariant ideals of $C^*(E)$ is finite, and there is a
finite chain $I_0=\{0\}\le I_1\le \cdots \le I_n=C^*(E) $ of
closed gauge-invariant ideals such that each quotient
$I_{i+1}/I_i$ is gauge-simple. We proceed by induction on $n$. If
$n=1$ we have the case in which $C^*(E)$ is gauge-simple, and thus
it is either purely infinite simple, or $AF$ or Morita-equivalent
to $C(\mathbb T)$; see \cite{BPRS}. In either case the result
follows. Note that in the purely infinite case, we use that
$V(C^*(E))=K_0(C^*(E))\setminus\{0\}=K_0(L)\setminus \{0\}=V(L)$.
Now assume that the result is true for graph $C^*$-algebras of
(gauge) length $n-1$ and let $A=C^*(E)$ be a graph $C^*$-algebra
of length $n$. Let $H$ be the saturated hereditary subset of $E^0$
corresponding to the ideal $I_1$. Note that $H$ is a minimal
saturated hereditary subset of $E^0$, and thus $H$ is cofinal. Set
$B=A/I_1$. By \cite[Theorem 4.1(b)]{BPRS}, we have $B\cong
C^*(F)$, where $F^0=E^0\setminus H$ and $F^1=\{e\in E^1\mid
r(e)\notin H \}$. Observe that by the induction hypothesis we know
that every projection in $B$ is equivalent to a finite orthogonal
sum of basic projections of the form $p_v$, where $v$ ranges in
$F^0=E^0\setminus H$. Let $\pi \colon A\to B $ denote the
canonical projection. Since $I_1$ is the closed ideal generated by
its projections, there is an embedding $V(A)/V(I_1)\to V(B)$. This
follows from \cite[Proposition 5.3(c)]{AF}, taking into account
that every closed ideal generated by projections is an almost
trace ideal. By induction hypothesis, $V(B)=V(C^*(F))$ is
generated as a monoid by $[p_v]$, for $v\in E^0\setminus H$, and
so the map $V(A)/V(I_1)\to V(B)$ is also surjective, so that
$V(B)\cong V(A)/V(I_1)$. In particular, $\pi (P)\sim \pi (Q)$ for
two projections $P,Q\in M_{\infty}(A)$, if and only if there are
projections $P',Q'\in M_{\infty}(I_1)$ such that $P\oplus P'\sim
Q\oplus Q'$ in $M_{\infty}(A)$.

\smallskip

We first deal with the case where $I_1$ has stable rank one, which
corresponds to the cases where $I_1$ is either AF or Morita
equivalent to $C(\mathbb T)$. Note that in this case either $H$
contains a sink $v$, or we have a simple loop without exit, in
which case we select $v$ as a vertex in this loop. Note that, by
the cofinality of $H$, any projection in $I_1$ is equivalent to a
projection of the form $k\cdot p_v$ for some $k\ge 0$. Now take
any projection $P$ in $M_{\infty}(A)$. Since $\pi (P)\sim
\pi(p_{v_1}\oplus \cdots \oplus p_{v_r})$ for some vertices
$v_1,\dots ,v_r$ in $E\setminus H$, there are $a,b\ge 0$ such that
$$P\oplus a\cdot p_v\sim p_{v_1}\oplus \cdots \oplus p_{v_r}\oplus b\cdot p_v.$$
Since the stable rank of $p_vAp_v$ is one, the projection $p_v$
cancels in direct sums \cite{rieffel}, and so, if $b\ge a$, we get
$$P\sim p_{v_1}\oplus \cdots \oplus p_{v_r}\oplus (b-a)p_v ,$$
so that $P$ is equivalent to a finite orthogonal sum of basic
projections. If $b<a$, then we have $P\oplus (a-b)p_v\sim
p_{v_1}\oplus \cdots \oplus p_{v_r} $. We claim that there is some
$1\le i\le r$ such that $v$ is in the tree of $v_i$. For, assume
to the contrary that $v\notin \bigcup _{i=1}^r T(v_i)$. We will
see that $v$ is not in the saturated hereditary subset of $E$
generated by $v_1,\dots ,v_r$. Note that the set $D=\bigcup
_{i=1}^r T(v_i)$ is hereditary, and that the saturated hereditary
subset of $E$ generated by $v_1,\dots ,v_r$ is $\ol{D}=\bigcup
_{j=0}^{\infty}\Lambda _i(D)$, see Section 5. Observe also that,
since $v$ is either a sink or belongs to a simple loop and $H$ is
cofinal, $v$ belongs to the tree of any vertex in $H$, whence
$H\cap D=\emptyset$. Let $v'$ be a vertex in $H$. If $v'\in
\Lambda _1(D)$ then $s^{-1}(v')\ne \emptyset$ and
$r(s^{-1}(v'))\subseteq D\cap H$. Since $H\cap D=\emptyset$, this
is impossible. So $\Lambda _1(D)\cap H=\emptyset$. Indeed, an easy
induction shows that $\Lambda _i(D)\cap H=\emptyset$ for all $i$,
and so $\ol{D}\cap H=\emptyset$. But being $p_v$ equivalent to a
subprojection of $p_{v_1}\oplus \cdots \oplus p_{v_r}$, the
projection $p_v$ belongs to the closed ideal of $A$ generated by
$p_{v_1}, \dots, p_{v_r}$, and so $v$ belongs to $\ol{D}$. This
contradiction shows that $v$ belongs to the tree of some $v_i$, as
claimed.

Now, the fact that $v$ belongs to the tree of $v_i$ implies that
there is a projection $Q$ which is a finite orthogonal sum of
basic projections such that $p_{v_i}\sim p_v\oplus Q$. Therefore
we get
$$P\oplus (a-b)p_v\sim (p_{v_1}\oplus \cdots p_{v_{i-1}}\oplus p_{v_{i+1}}\oplus \cdots
\oplus p_{v_r}\oplus Q)\oplus p_v . $$ Since $p_v$ can be
cancelled in direct sums, we get
$$P\oplus (a-b-1)p_v\sim (p_{v_1}\oplus \cdots p_{v_{i-1}}\oplus p_{v_{i+1}}\oplus \cdots
\oplus p_{v_r}\oplus Q), $$ and so, using induction, we obtain
that $P$ is equivalent to a finite orthogonal sum of basic
projections.

\smallskip

Finally we consider the case where $I_1$ is a purely infinite
simple $C^*$-algebra. Recall that in this case  $I_1$ has real
rank zero \cite{BP}, and that $V(I_1)\setminus \{ 0\}$ is a group.
So there is a nonzero projection $e$ in $I_1$ such that for every
nonzero projection $p$ in $M_{\infty}(I_1)$ there exists a nonzero
projection $q\in I_1$ such that $p\oplus q\sim e$. Let $P$ be a
nonzero projection in $M_k(A)$, for some $k\ge 1$, and denote by
$I$ the closed ideal of $A$ generated by (the entries of) $P$. If
$I\cdot I_1=0$, then $I\cong (I+I_1)/I_1$, so that $I$ is a closed
ideal in the quotient $C^*$-algebra $B=A/I_1$. It follows then by
our assumption on $B$ that $P$ is equivalent to a finite
orthogonal sum of basic projections. Assume now that $I\cdot
I_1\ne 0$. Then there is a nonzero column $C=(a_1,a_2,\dots
,a_k)^t\in A^k$ such that $C=PCe$. Consider the positive element
$c=C^*C$ , which belongs to $eAe$. Since $e\in I_1$ and $I_1$ has
real rank zero, the $C^*$-algebra $eAe$ has also real rank zero,
so that we can find a nonzero projection $p\in cAc$. Take $x\in A$
such that $p=cxc$. By using standard tricks (see e.g.
\cite{rordam}), we can now produce a projection $P'\le P$ such
that $p\sim P'$. Namely, consider the idempotent $F=CpC^*CxC^*$ in
$PM_k(A)P$. Then $p$ and $F$ are equivalent as idempotents, and
$F$ is equivalent to some projection $P'$ in $PM_k(A)P$; see
\cite[Exercise 3.11(i)]{rordam}. Since $p$ and $P'$ are equivalent
as idempotents, they are also Murray-von Neumann equivalent, see
\cite[Exercise 3.11(ii)]{rordam}, as desired. We have proved that
there is a nonzero projection $p$ in $I_1$ such that $p$ is
equivalent to a subprojection of $P$. Since $I_1$ is purely
infinite simple, every projection in $I_1$ is equivalent to a
subprojection of $p$, and so every projection in $I_1$ is
equivalent to a subprojection of $P$.

\smallskip

Now we are ready to conclude the proof. There is a projection $q$
in $I_1$ such that $P\oplus q$ is equivalent to a finite
orthogonal sum of basic projections. Let $q'$ be a nonzero
projection in $I_1$ such that $q\oplus q'\sim e$, and observe that
$$P\oplus e\sim (P\oplus q)\oplus q' ,$$
so that $P\oplus e$ is also a finite orthogonal sum of basic
projections. By the above argument, there is a projection $e'$
such that $e'\le P$ and $e\sim e'$. Write $P=e'+P'$. Then we have
$$P\oplus e\sim P'\oplus e'\oplus e\sim P'\oplus e\oplus e\sim P'\oplus e\sim P.$$
It follows that $P\sim P\oplus e$ and so $P$ is equivalent to a
finite orthogonal sum of basic projections.
\end{proof}

\begin{corollary}
\label{stableweak} Let $E$ be a row-finite graph. Then the monoid
$V(C^*(E))$ is a refinement monoid and $C^*(E)$ has stable weak
cancellation. Moreover, $V(C^*(E))$ is an unperforated monoid and
$K_0(C^*(E))$ is an unperforated group.
\end{corollary}

\begin{proof}
By Theorem \ref{C*algebras}, $V(C^*(E))\cong M_E$, and so
$V(C^*(E))$ is a refinement monoid by Proposition
\ref{refinement}.

It follows from Theorem \ref{separativity} that $V(C^*(E))$ is a
separative monoid. By Proposition \ref{swc}, this is equivalent to
saying that $C^*(E)$ has stable weak cancellation. The statements
about unperforation follow from Proposition \ref{unperforation}.
\end{proof}

\section*{Acknowledgments}

Part of this work was done during a visit of the second author to
the the Departament de Matem\`atiques de l'Universitat Aut\`onoma
de Barcelona, and during a visit of the third author to the Centre
de Recerca Matem\`atica (UAB). The second and third authors want
to thank both host centers for their warm hospitality. The first
author thanks Mikael R\o rdam for valuable discussions on the
topic of Section 7. The authors thank Ken Goodearl for interesting
discussions at an early stage in the preparation of this work, and
Fred Wehrung for several useful comments and for pointing out to
them the papers \cite{Brook} and \cite{Wehrung}.


\begin{thebibliography}{99}

\bibitem{AA} \textsc{G. Abrams, G. Aranda Pino},
The Leavitt path algebra of a graph, \emph{J. Algebra} {\bf 293}
(2005), 319--334.

\bibitem{AGOP} \textsc{P. Ara, K. R. Goodearl, K. C. O'Meara, E. Pardo},
Separative cancellation for projective modules over exchange
rings, \emph{Israel J. Math.} {\bf 105} (1998), 105--137.

\bibitem{AGOR}\textsc{P. Ara, K. R. Goodearl, K. C. O'Meara, R.
Raphael}, $K_1$ of separative exchange rings and $C^*$-algebras
with real rank zero, \emph{Pacific J. Math.} {\bf 195} (2000),
261--275.

\bibitem{AF} \textsc{P. Ara, A. Facchini}, Direct sum
decompositions of modules, almost trace ideals, and pullbacks of
monoids, \emph{Forum Math.} {\bf 18} (2006), 365--389.


\bibitem{BPRS} \textsc{T. Bates, D. Pask, I. Raeburn,
W. Szyma\' nski}, The $C^*$-algebras of row-finite graphs,
\emph{New York J. Math.} \textbf{6} (2000), 307--324.

\bibitem{Bergman} \textsc{G.M. Bergman}, Coproducts and some
universal ring constructions, \emph{Trans. Amer. Math. Soc.}
\textbf{200} (1974), 33--88.

\bibitem{RatC*} {\sc B. Blackadar}, Rational C*-algebras and
nonstable K-Theory, \emph{Rocky Mountain J. Math.} {\bf 20(2)}
(1990), 285--316.

\bibitem{Black} {\sc B. Blackadar}, ``K-Theory for Operator
Algebras'', Second Edition, M.S.R.I. Publications, vol. 5,
Cambridge Univ. Press, Cambridge, 1998.

\bibitem{Brook}{\sc G. Brookfield}, Cancellation in primely
generated refinement monoids, \emph{Algebra Universalis}, {\bf 46}
(2001), 343--371.

\bibitem{brown} \textsc{L. G. Brown}, Homotopy of projections in $C\sp *$-algebras
of stable rank one, in \emph{Recent advances in operator algebras}
(Orl\'{e}ans, 1992). Ast\'{e}risque No. 232 (1995), 115--120.

\bibitem{BP} \textsc{L. G. Brown, G. K. Pedersen}, $C^{*}$-algebras of
real rank zero, \emph{J. Funct. Anal.} \textbf{99} (1991),
131--149.

\bibitem{BPns} \textsc{L. G. Brown, G. K. Pedersen}, Non-stable $K$-theory
and extremally rich $C^*$-algebras, \emph{Preprint}.

\bibitem{CK} \textsc{J. Cuntz, W. Krieger}, A class of
C*-algebras and topological Markov chains, \emph{Inventiones
Math.} \textbf{56} (1980), 251--268.

\bibitem{DHS} \textsc{K. Deicke, J.H. Hong, W. Szyma\' nski},
Stable rank of graph algebras. Type I graph algebras and their
limits, \emph{Indiana Univ. Math. J.}  \textbf{52(4)} (2003),
963--979.


\bibitem{FHK} \textsc{A. Facchini, F. Halter-Koch}, Projective modules and
divisor homomorphisms, \emph{J. Algebra Appl.} {\bf 2} (2003),
435--449.

\bibitem{vnrr} \textsc{K. R. Goodearl},
``Von Neumann Regular Rings'', Second edition, Krieger Publishing
Co., Inc., Malabar, FL, 1991.



\bibitem{JP} \textsc{J.A. Jeong, G.H. Park}, Graph C*-algebras
with real rank zero, \emph{J. Funct. Anal.} \textbf{188} (2002),
216--226.


\bibitem{KPR} \textsc{A. Kumjian, D. Pask, I. Raeburn}, Cuntz-Krieger algebras
of directed graphs, \emph{Pacific J. Math.} \textbf{184} (1998),
161--174.

\bibitem{Leavitt} \textsc{W. G. Leavitt}, The module type of a
ring, \emph{Trans. Amer. Math. Soc.} \textbf{42} (1962), 113--130.

\bibitem{MM} \textsc{P. Menal, J. Moncasi}, Lifting units in
self-injective rings and an index theory for Rickart
$C^*$-algebras, \emph{Pacific J. Math.} {\bf 126} (1987),
295--329.

\bibitem{Perera}\textsc{F. Perera}, Lifting units modulo exchange
ideals and $C^*$-algebras with real rank zero, \emph{J. reine
angew. Math.} {\bf 522} (2000), 51--62.

\bibitem{RS} \textsc{I. Raeburn, W. Szyma\' nski},
Cuntz-Krieger algebras of infinite graphs and matrices,
\emph{Trans. Amer. Math. Soc.} \textbf{356} (2004), 39--59.

\bibitem{rieffel} \textsc{M.A. Rieffel}, Dimension and stable rank in the $K$-theory of
$C\sp{*}$-algebras, \emph{Proc. London Math. Soc.} \textbf{46}
(1983), 301--333.

\bibitem{rordam}\textsc{M. R\o rdam, F. Larsen, N.J. Laustsen},
``An Introduction to $K$-Theory for $C^*$-Algebras", Cambridge
University Press, LMS Student Texts 49, 2000.

\bibitem{Ros} \textsc{J. Rosenberg}, ``Algebraic $K$-Theory and Its
Applications", Springer-Verlag, GTM 147, 1994.

\bibitem{Wat} \textsc{Y. Watatani}, Graph theory for
$C^*$-algebras, in {\em Operator algebras and their applications}
(R. V. Kadison, ed.), Proc. Sympos. Pure Math., vol {\bf 38} Part
I, Amer. Math. Soc., Providence, 1982, pp. 195--197.

\bibitem{Wehrung}\textsc{F. Wehrung}, The Dimension Monoid of a
Lattice, \emph{Algebra Universalis} \textbf{40} (1998), 247--411.



\bibitem{zhang} \textsc{S. Zhang}, Diagonalizing projections
in multiplier algebras and in matrices over a $C^*$-algebra,
\emph{Pacific J. Math.} \textbf{54} (1990), 181--200.
\end{thebibliography}
\end{document}